\documentclass[a4paper,12pt]{amsart}

\usepackage{amssymb}
\usepackage{bbm}
\usepackage{mathrsfs}	
\usepackage[all]{xy}
\usepackage[T5,T1]{fontenc}

\normalsize

\newtheoremstyle{bracket}{1ex}{2ex}{\rm}{}{\bfseries}{}{0.8em}{\thmnumber{(#2)}}
\newtheoremstyle{thm}{1ex}{2ex}{\itshape}{}{\bfseries}{}{0.9em}{\thmnumber{(#2)}\thmname{ #1}\thmnote{ (#3)}}

\theoremstyle{bracket}
\newtheorem{no}{}

\theoremstyle{thm}
\newtheorem{prop}[no]{Proposition}
\newtheorem{theorem}[no]{Theorem}
\newtheorem{cor}[no]{Corollary}

\DeclareMathOperator{\spec}{Spec}
\DeclareMathOperator{\pic}{Pic}

\newcommand{\dfgl}{\mathrel{\mathop:}=}

\newcommand{\ftnt}{\footnotetext{This text is a slight update of the author's contribution to the proceedings of the 7th Japan-Vietnam Joint Seminar on Commutative Algebra, which itself was partially an update of parts of the author's contribution to the proceedings of the 6th Japan-Vietnam Joint Seminar on Commutative Algebra. The author was supported by the Swiss National Science Foundation.}}

\begin{document}

\title[On quasicoherent sheaves on toric schemes]{On quasicoherent sheaves\\on toric schemes\ftnt}
\author{Fred Rohrer}
\address{Universit\"at T\"ubingen, Fachbereich Mathematik, Auf der Morgenstelle 10, 72076 T\"ubingen, Germany}
\email{fredrohrer0@gmail.com}
\dedicatory{{\fontencoding{T5}\selectfont D\`anh ri\ecircumflex ng cho t\'\acircumflex t c\h a ng\uhorn\`\ohorn i b\d an Vi\d\ecircumflex t c\h ua t\ocircumflex i}}

\begin{abstract}
A correspondence between quasicoherent sheaves on toric schemes and graded modules over some homogeneous coordinate ring is presented, and the behaviour of several finiteness properties under this correspondence is investigated.
\end{abstract}

\maketitle


\section*{Introduction}

Toric varieties are probably well-known to anyone working in algebraic geometry. A lot is known about them, and a lot still gets published about them. Their theory was generalised in several directions, and this often also led to a better understanding of classical toric varieties. The most natural generalisation seems to be the \textit{study of toric varieties from a scheme-theoretical point of view.} It is clear that to do this, one has to be able to make arbitrary base changes. Hence, instead of considering toric varieties over an algebraically closed field (or, as often done, over the field of complex numbers $\mathbbm{C}$), one needs to study \textit{toric schemes,} i.e., ``toric varieties over arbitrary base rings''. With this motivation it is not hard to see that our toric schemes have to be defined as being constructed from a fan in the same way as the classical toric varieties. Special cases of this generalisation were mentioned briefly in \cite[\textsection 4]{dem} (for regular fans and mainly over the ring of integers). Of course, one could also consider the equivalent description of toric varieties over $\mathbbm{C}$ in terms of torus operations, and try to replace $\mathbbm{C}$ by more general base rings. The case of discrete valuation rings was treated in \cite[IV.3]{kkms}, and more general valuation rings are considered in recent work by Gubler \cite[Chapters 6--7]{gubler}. However, this approach leads to a different class of schemes and hence differs heavily from ours, as it requires for example the consideration of a bigger or otherwise different class of fans.

Besides yielding a better understanding of the geometry of toric varieties, there are concrete applications of the above generalisation, as the following remark shows. Let $X$ be the toric variety over an algebraically closed field $K$ associated with a fan $\Sigma$. A fundamental question in algebraic geometry is then whether the Hilbert functor ${\rm Hilb}_{X/K}$ of $X$ over $K$ is representable, i.e., whether the Hilbert scheme of $X$ exists (cf.~\cite{fga}). If $X$ is projective then this is indeed the case and follows from Grothendieck's more general result \cite[Th\'eor\`eme 3.1]{fga}. However, toric varieties are not necessarily projective (nor quasiprojective), and in general it is not known whether their Hilbert schemes exist. Studying ${\rm Hilb}_{X/K}$ amounts to studying quasicoherent sheaves on the base change $X\otimes_KR$ for every $K$-algebra $R$, and it turns out that $X\otimes_KR$ is precisely the toric scheme over $R$ associated with $\Sigma$. Hence, in order to study Hilbert functors of toric varieties \textit{it is necessary to study toric schemes over more general bases than just over algebraically closed fields.}

\smallskip

Here we give an overview of some results about quasicoherent sheaves on toric schemes. They generalise results by Cox \cite{cox} and Musta\c{t}\v{a} \cite{mus1}, and they extend some results presented at the 6th Japan-Vietnam Joint Seminar on Commutative Algebra \cite{ts0}. Furthermore, they are mostly inspired by well-known analogous results about projective schemes as given in \cite[Chapitre II]{ega}. The corresponding proofs are still unpublished$^0$\footnotetext{Not anymore -- see F.\,Rohrer, {\it Quasicoherent sheaves on toric schemes.} Expo.\,Math.\,32 (2014) 33--78.}.


\section{Toric schemes}

We start by briefly describing the construction of toric schemes from fans. For unexplained terminology and notation from polyhedral geometry we refer the reader to \cite{ts0}, \cite[Chapter II]{diss} (available from the author's homepage), or \cite{gts}.

\smallskip

\noindent\textit{$\bullet$\quad From now on let $V$ be an $\mathbbm{R}$-vector space of finite dimension $n$, let $N$ be a $\mathbbm{Z}$-structure on $V$ (i.e., a subgroup of rank $n$ of the additive group underlying $V$ with $\langle N\rangle_{\mathbbm{R}}=V$), let $M\dfgl N^*$ denote the dual of $N$, which is a $\mathbbm{Z}$-structure on the dual $V^*$ of $V$, let $\Sigma$ be an $N$-fan in $V$, and let $R$ be a ring\/\footnote{Rings, groups and monoids are understood to be commutative, and algebras are understood to be commutative, unital and associative.}.}

\smallskip

If $\sigma\in\Sigma$ then $\sigma^{\vee}\cap M$ is a torsionfree, cancellable, finitely generated submonoid of $M$, and if $\tau$ is a face of $\sigma$ then $\sigma^{\vee}\cap M$ is a submonoid of $\tau^{\vee}\cap M$. Taking spectra of algebras of monoids over $R$ and setting $X_{\sigma}(R)\dfgl\spec(R[\sigma^{\vee}\cap M])$ for $\sigma\in\Sigma$ we get an inductive system $$\bigl(X_{\sigma}(R)\rightarrow\spec(R)\bigr)_{\sigma\in\Sigma}$$ of $R$-schemes over $\Sigma$. Its inductive limit exists and is an $R$-scheme denoted by $$X_{\Sigma}(R)\rightarrow\spec(R)$$ and called \textit{the toric scheme over $R$ associated with $\Sigma$ (and $N$).} It can be understood as obtained by gluing $(X_{\sigma}(R))_{\sigma\in\Sigma}$ along $(X_{\sigma\cap\tau}(R))_{(\sigma,\tau)\in\Sigma^2}$.

The above construction of toric schemes gives rise to a contravariant functor $X_{\Sigma}$ from the category of rings to the category of schemes together with a morphism $X_{\Sigma}\rightarrow\spec$. This functor is compatible with base change in the sense that there is a canonical isomorphism $$X_{\Sigma}(\bullet)\cong X_{\Sigma}(R)\otimes_R\bullet$$ of contravariant functors from the category of $R$-algebras to the category of $R$-schemes.

\smallskip

The first important question is now of course how the base ring affects the geometry of a toric scheme. It turns out that all toric schemes share some basic properties and thus are not too ugly. For example, the $R$-scheme $X_{\Sigma}(R)\rightarrow\spec(R)$ is separated, quasicompact, flat, and of finite presentation, and it is faithfully flat if and only if $\Sigma\neq\emptyset$ or $R=0$. However, a lot of other basic properties are respected and reflected by $X_{\Sigma}$. For example, the scheme $X_{\Sigma}(R)$ is reduced, connected, normal, or Noetherian if and only if $R$ is so or $\Sigma=\emptyset$, and it is irreducible or integral if and only if $R$ is so and $\Sigma\neq\emptyset$. Finally, some properties depend also on the fan. For example, the $R$-scheme $X_{\Sigma}(R)\rightarrow\spec(R)$ is proper if and only if $\Sigma$ is complete, or $\Sigma=\emptyset$, or $R=0$. (More results of this type can be found in \cite{ts0} and \cite{gts}.)

The above examples show in particular that on general toric schemes \textit{no satisfying theory of Weil divisors is available.} Since a lot of results about toric varieties were proved by heavy use of Weil divisor techniques (see e.g.~\cite{cox}, \cite{ful}), one has to come up with new proofs in order to generalise these results to toric schemes.


\section{Sheaves on toric schemes}

Generalising work of Cox \cite{cox} and Musta\c{t}\v{a} \cite{mus1} we introduce a notion of Cox ring (not to be confused with the one introduced in \cite{hukeel}) and describe quasicoherent modules on toric schemes in terms of graded modules over these rings. In order to do so we need to define some objects encoding the combinatorics of the fan $\Sigma$.

\smallskip

We denote by $\Sigma_1$ the set of $1$-dimensional cones in $\Sigma$, and for $\rho\in\Sigma_1$ by $\rho_N$ its unique minimal $N$-generator. There is an exact sequence of groups $$M\overset{c}\longrightarrow\mathbbm{Z}^{\Sigma_1}\overset{a}\longrightarrow A\longrightarrow 0,$$ where $c(u)\dfgl(u(\rho_N))_{\rho\in\Sigma_1}$ for $u\in M$ and where $a$ is defined as the cokernel of $c$. Note that $c$ is a monomorphism if and only if $\Sigma$ is \textit{full,} i.e., $\langle\bigcup\Sigma\rangle_{\mathbbm{R}}=V$. We denote by $(\delta_{\rho})_{\rho\in\Sigma_1}$ the canonical basis of $\mathbbm{Z}^{\Sigma_1}$ and set $\alpha_{\rho}\dfgl a(\delta_{\rho})$ for $\rho\in\Sigma_1$.

\smallskip

\noindent\textit{$\bullet$\quad From now on let $B\subseteq A$ be a subgroup.}

\smallskip

We denote by $S$ the polynomial algebra $R[(Z_{\rho})_{\rho\in\Sigma_1}]$ in indeterminates $(Z_{\rho})_{\rho\in\Sigma_1}$ over $R$, furnished with the $A$-graduation induced by $a$, i.e., such that $\deg(Z_{\rho})=\alpha_{\rho}$ for $\rho\in\Sigma_1$. For $\sigma\in\Sigma$ we set $\widehat{Z}_{\sigma}\dfgl\prod_{\rho\in\Sigma_1\setminus\sigma_1}Z_{\rho}\in S$ (where $\sigma_1$ denotes the set of $1$-dimensional faces of $\sigma$). Finally we define a graded ideal $I\dfgl\langle\widehat{Z}_{\sigma}\mid\sigma\in\Sigma\rangle_S$.

The $B$-graded $R$-algebra $S_B\dfgl\bigoplus_{\alpha\in B}S_{\alpha}$ obtained from $S$ by degree restriction to $B$ is called \textit{the $B$-restricted Cox ring over $R$ associated with $\Sigma_1$ (and $N$),} and its graded ideal $I_B\dfgl I\cap S_B$ is called \textit{the $B$-restricted irrelevant ideal over $R$ associated with $\Sigma$ (and $N$);} it is generated by finitely many monomials.

To proceed we need to ``invert the monomials $\widehat{Z}_{\sigma}$ in the $B$-restricted Cox ring'', and hence we have to assure that some powers of these monomials lie in $S_B$. This amounts to requiring $B$ to be \textit{big,} i.e., to have finite index in $A$.

\smallskip

\noindent\textit{$\bullet$\quad From now on suppose that $B$ is big, so that there exists $m\in\mathbbm{N}_0$ with $\widehat{Z}_{\sigma}^m\in S_B$ for every $\sigma\in\Sigma$.}

\smallskip

For $\sigma\in\Sigma$ the $B$-graded ring of fractions $(S_B)_{\widehat{Z}_{\sigma}^m}$ is independent of the choice of $m$, and it is denoted by $S_{B,\sigma}$. Its component of degree $0$ is independent of the choice of $B$ and is denoted by $S_{(\sigma)}$. Moreover, for a face $\tau$ of $\sigma$ there is a canonical morphism of rings $S_{(\sigma)}\rightarrow S_{(\tau)}$, which is independent of $m$ and $B$. Taking spectra and setting $Y_{(\sigma)}(R)\dfgl\spec(S_{(\sigma)})$ for $\sigma\in\Sigma$ we obtain an inductive system $$\bigl(Y_{\sigma}(R)\rightarrow\spec(R)\bigr)_{\sigma\in\Sigma}$$ of $R$-schemes over $\Sigma$. Its inductive limit exists and is an $R$-scheme denoted by $$Y_{\Sigma}(R)\rightarrow\spec(R)$$ and called \textit{the Cox scheme over $R$ associated with $\Sigma$ (and $N$).} It can be understood as obtained by gluing $(Y_{\sigma}(R))_{\sigma\in\Sigma}$ along $(Y_{\sigma\cap\tau}(R))_{(\sigma,\tau)\in\Sigma^2}$.

The above construction of Cox schemes gives rise to a contravariant functor $Y_{\Sigma}$ from the category of rings to the category of schemes together with a morphism $Y_{\Sigma}\rightarrow\spec$. This functor is compatible with base change in the sense that there is a canonical isomorphism $$Y_{\Sigma}(\bullet)\cong Y_{\Sigma}(R)\otimes_R\bullet$$ of contravariant functors from the category of $R$-algebras to the category of $R$-schemes.

\smallskip

Cox schemes are closely related to toric schemes as follows. The morphism of groups $c:M\rightarrow\mathbbm{Z}^{\Sigma_1}$ induces for $\sigma\in\Sigma$ a morphism of rings $R[\sigma^{\vee}\cap M]\rightarrow S_{(\sigma)}$, and these morphisms induce a canonical morphism of contravariant functors $\gamma:Y_{\Sigma}\rightarrow X_{\Sigma}$. It turns out that this morphism is an isomorphism if and only if $\Sigma$ is full.

Therefore, using a straightforward (but necessarily noncanonical) procedure to consider a toric scheme associated with a nonfull fan as a toric scheme associated with a full fan it is sufficient to study henceforth Cox schemes instead of toric schemes. (Note that this reduction demands a base change and is in general \textit{not available for toric varieties.})

\smallskip

Now we are ready to explain how $B$-graded $S_B$-modules give rise to quasicoherent sheaves on $Y\dfgl Y_{\Sigma}(R)$. We denote by ${\sf GrMod}^B(S_B)$ and ${\sf QCMod}(\mathscr{O}_Y)$ the categories of $B$-graded $S_B$-modules and of quasicoherent $\mathscr{O}_Y$-modules, respectively. Moreover, for a $B$-graded $S_B$-module $F$ we denote by $F_{(\sigma)}$ the component of degree $0$ of the $B$-graded module of fractions $F_{\widehat{Z}_{\sigma}^m}=F\otimes_{S_B}(S_B)_{\widehat{Z}_{\sigma}^m}$, and for an $S_{(\sigma)}$-module $G$ we denote by $\widetilde{G}$ the $\mathscr{O}_{Y_{\sigma}(R)}$-module associated with $G$. There exists a unique functor $$\mathscr{S}_B:{\sf GrMod}^B(S_B)\rightarrow{\sf QCMod}(\mathscr{O}_Y)$$ with $\mathscr{S}_B(F)\!\upharpoonright_{Y_{\sigma}(R)}=\widetilde{F_{(\sigma)}}$ for every $\sigma\in\Sigma$ and every $B$-graded $S_B$-module $F$. Since $\mathscr{S}_B$ coincides locally with the well-known equivalence between modules and quasicoherent sheaves on affine schemes it is exact and commutes with inductive limits.

Denoting by $\bullet(\alpha)$ the functor of shifting degrees by $\alpha$, we can construct a functor $$\Gamma_*^B(\bullet)\dfgl\bigoplus_{\alpha\in B}\Gamma\bigl(Y,\bigl(\bullet\otimes_{\mathscr{O}_Y}\mathscr{S}_B(S_B(\alpha))\bigr)\bigr)$$ for $\mathscr{S}_B$, called \textit{the first total functor of sections associated with $\Sigma$ and $B$ over $R$,} and show that there is an isomorphism of functors $$\beta_B:\mathscr{S}_B\circ\Gamma_*^B\rightarrow{\rm Id}_{{\sf GrMod}^B(S_B)}.$$ Thus, we get the following generalisation of \cite[Theorem 1.1]{mus1}, itself a generalisation of \cite[Theorem 3.2]{cox}.

\begin{theorem}\label{surj}
The functor $\mathscr{S}_B:{\sf GrMod}^B(S_B)\rightarrow{\sf QCMod}(\mathscr{O}_Y)$ is essentially surjective.
\end{theorem}

Keeping in mind that the functor $\mathscr{S}_B$ is exact (and also keeping in mind the analogous results in the case of projective schemes) we look now at the following situation. Let $F$ be a $B$-graded $S_B$-module. A graded sub-$S_B$-module $G\subseteq F$ is said to be $I_B$-saturated (in $F$) if $G=\bigcup_{m\in\mathbbm{N}}(G:_FI_B^m)$, where the right hand side can be shown to be the smallest graded sub-$S_B$-module of $F$ containing $G$ that is $I_B$-saturated. We denote by $\widetilde{\mathbbm{J}}_F$, $\mathbbm{J}_F$, and $\mathbbm{J}_F^{{\rm sat}}$ the sets of quasicoherent sub-$\mathscr{O}_Y$-modules of $\mathscr{S}_B(F)$, of graded sub-$S_B$-modules of $F$, and of $I_B$-saturated graded sub-$S_B$-modules of $F$, respectively. Then, the functor $\mathscr{S}_B$ induces a map $\Xi_F:\mathbbm{J}_F\rightarrow\widetilde{\mathbbm{J}}_F$, and this map restricts to a map $\Xi_F^{{\rm sat}}:\mathbbm{J}_F^{{\rm sat}}\rightarrow\widetilde{\mathbbm{J}}_F$.

The question is now whether the maps $\Xi_F$ or $\Xi_F^{{\rm sat}}$ are surjective or injective. In case $F=S_B$, i.e., when we consider only ideals, this was already settled in \cite{ts0}. There, the result about injectivity required the big subgroup $B$ not to be ``too big''. More precisely, $B$ is called \textit{small (with respect to $\Sigma$)} if it is contained in the Picard group $\pic(\Sigma)$ of $\Sigma$, as introduced in \cite[V.5]{ewald} (see also \cite{ts0}). Then, we gave the following result.

\begin{prop}\label{ideal}
a) The map $\Xi_{S_B}:\mathbbm{J}_{S_B}\rightarrow\widetilde{\mathbbm{J}}_{S_B}$ is surjective.

b) The map $\Xi_{S_B}^{{\rm sat}}:\mathbbm{J}_{S_B}^{{\rm sat}}\rightarrow\widetilde{\mathbbm{J}}_{S_B}$ is surjective, and if $B$ is small then it is bijective.
\end{prop}

To get a generalisation for arbitrary $F$ we need a further ingredient, namely the second total functor of sections. (In \cite{ts0} this was introduced and used in the context of cohomology to get the toric Serre-Grothendieck correspondence.) We define a functor $$\Gamma_{**}^B(\bullet):{\sf GrMod}^B(S_B)\rightarrow{\sf GrMod}^B(S_B),$$ called \textit{the second total functor of sections associated with $\Sigma$ and $B$ over $R$,} by setting $$\Gamma_{**}^B(\bullet)\dfgl\bigoplus_{\alpha\in B}\Gamma\bigl(Y,\mathscr{S}_B(\bullet(\alpha))\bigr).$$ (Despite its name it is defined on the category ${\sf GrMod}^B(S_B)$, but by (\ref{surj}) this is merely a technical point.) The reason that there are two (in general different) total functors of sections is that the canonical morphism $$\mathscr{S}_B(\bullet)\otimes_{\mathscr{O}_{Y_{\Sigma}(R)}}\mathscr{S}_B(S_B(\alpha))\rightarrow\mathscr{S}_B(\bullet(\alpha))$$ is not necessarily an isomorphism, and thus the induced morphism of functors $$\delta_B:\Gamma_*^B\circ\mathscr{S}_B\rightarrow\Gamma_{**}^B$$ is not necessarily an isomorphism. However, after composing it with $\mathscr{S}_B$ it turns into an isomorphism. More precisely, there is a canonical morphism of functors $$\eta_B:{\rm Id}_{{\sf GrMod}^B(S_B)}\rightarrow\Gamma_{**}^B,$$ and we have the following result.

\begin{prop}
The diagram of functors $$\xymatrix{
&\mathscr{S}_B\circ\Gamma_*^B\circ\mathscr{S}_B\ar[ld]_{\beta_B\circ\mathscr{S}_B}\ar[rd]^{\mathscr{S}_B\circ\delta_B}&\\
\mathscr{S}_B\ar[rr]^{\mathscr{S}_B\circ\eta_B}&&\mathscr{S}_B\circ\Gamma_{**}^B
}$$ commutes, and all occuring morphisms are isomorphisms.
\end{prop}

This result is the key observation for the following generalisation of (\ref{ideal}).

\begin{theorem}\label{bij}
Let $F$ be a $B$-graded $S_B$-module.

a) The map $\Xi_F:\mathbbm{J}_F\rightarrow\widetilde{\mathbbm{J}}_F$ is surjective.

b) The map $\Xi_F^{{\rm sat}}:\mathbbm{J}_F^{{\rm sat}}\rightarrow\widetilde{\mathbbm{J}}_F$ is surjective, and if $B$ is small then it is bijective.

c) If $F$ is an $I_B$-torsion module then it holds $\mathscr{S}_B(F)=0$, and if $B$ is small then the converse is true.
\end{theorem}

Keeping in mind that $\pic(\Sigma)$ is big if and only if $\Sigma$ is simplicial and that $\pic(\Sigma)\cong\pic(X_{\Sigma}(\mathbbm{C}))$ (\cite[Theorem VII.2.15]{ewald}), we get back \cite[Corollary 3.9]{cox} as a special case.


\section{Finiteness conditions}

Even when working in great generality, it is sometimes convenient to know something about finiteness conditions. Motivated by the fact that a homogeneous coordinate ring of a projective scheme can always be chosen such that it is of finite type over the base ring, we start by giving a condition such that the restricted Cox ring $S_B$ is of finite type over $R$. Recall that a $B$-graded $R$-algebra $T$ is said to be positively graded if $T_{\alpha}\neq 0$ and $T_{-\alpha}\neq 0$ for $\alpha\in B$ implies $\alpha=0$.

\begin{prop}
If the cone and the vector space generated by $\Sigma$ coincide, then $S_B$ is a positively $B$-graded $R$-algebra of finite type.
\end{prop}

Next we ask how certain finiteness conditions behave under the functor $\mathscr{S}_B$. Making use of (\ref{bij}) and quasicompacity of $Y$ we first get the following result.

\begin{prop}
If $F$ is a $B$-graded $S_B$-module and $\mathscr{G}\in\widetilde{\mathbbm{J}}_F$ is of finite type, then there exists $G\in\mathbbm{J}_F$ of finite type with $\mathscr{S}_B(G)=\mathscr{G}$.
\end{prop}

If we want to know whether $\mathscr{S}_B$ preserves some finiteness conditions, it turns out that we need to make some hypothesis on the subgroup $B$. More precisely, we need the graded rings of fractions $S_{B,\sigma}$ for $\sigma\in\Sigma$ to be strongly graded, and this is fulfilled if $B$ is small. Note that since $B$ is big this implies that $\Sigma$ is simplicial. The finiteness conditions that we consider here are being of finite type or of finite presentation, and being pseudocoherent or coherent. Recall that a graded module or a sheaf of modules is called pseudocoherent if its subobjects of finite type are of finite presentation, and coherent if it is pseudocoherent and of finite type. (For sheaves this definition coincides with the usual one found in \cite[Chapitre 0]{ega}.)

\begin{prop}\label{fin}
Suppose that $B$ is small, and let $F$ be a $B$-graded $S_B$-module. If $F$ is of finite type, of finite presentation, pseudocoherent, or coherent, then so is $\mathscr{S}_B(F)$.
\end{prop}

As an application of (\ref{fin}) we get a criterion for the structure sheaf $\mathscr{O}_Y$ to be coherent. Recall that a ring $T$ is called stably coherent if polynomial algebras in finitely many indeterminates over $T$ are coherent. Stable coherence is strictly weaker than coherence, but it is shared for example by absolutely flat rings, valuation rings, or semihereditary rings.

\begin{cor}
If $\Sigma$ is $N$-regular and $R$ is stably coherent, then $\mathscr{O}_Y$ is coherent.
\end{cor}

(The regularity hypothesis comes from the fact that it is equivalent to $\pic(\Sigma)=A$, so that choosing $B=\pic(\Sigma)=A$ it follows that $S_B$ is a polynomial algebra over $R$.)


\end{document}